\documentclass[11pt,Color]{article}
\usepackage{amssymb}
\usepackage{amsmath}
\usepackage{color}
\usepackage{indentfirst}

\newcommand{\mc}[1]{\mathcal{#1}}
\newcommand{\mf}[1]{\mathfrak{#1}}

\input xypic.sty

\def\qed{\hfill {\large ${\sqcup\!\!\!\!\sqcap}$}}

\newenvironment{demo}{{\bf Proof }}
{\qed \\}

\newcommand{\re}{\mathbb R}

\newcommand{\eme}{{\rm M}^{n+1}_\lambda}
\newcommand{\flecha}{\longrightarrow}
\newcommand{\<}{\left<}

\renewcommand{\>}{\right>}
\newcommand{\eps}{\ensuremath{\varepsilon}}
\def\bal{\begin{align}}
\def\eal{\end{align}}

\newcommand{\bde}{\begin{defi}}
\newcommand{\ede}{\end{defi}}

\renewcommand{\r}{\tilde{r}_t}

\numberwithin{equation}{section}
\def\be{\begin{equation}}
\def\ee{\end{equation}}
\def\mad{{\rm maxd}}
\def\H{{\mathcal H}}
\def\oH{{\overline H}}
\def\oR{{\overline R}}
\def\oRic{{\overline Ric}}
\def\ta{{\rm ta_\lambda}}
\def\a{\alpha}
\def\s{{\rm s_\lambda}}
\def\c{{\rm c_\lambda}}
\def\co{{\rm co_\lambda}}
\def\ral{\sqrt{|\lambda|}}
\def\tr{{\rm tr}}
\def\det{{\rm det}}
\def\dist{{\rm dist}}

\def\grad{{\rm grad}}
\def\grada{\overline{\rm{grad}}\ }
\def\nablaa{\overline{\nabla}}
\def\Deltaa{\overline{\Delta}}
\def\vle{{\rm vol}}
\def\parcial#1#2{\frac{\partial #1}{\partial#2}}

\def\Daparcial#1#2{\frac{\nablaa #1}{\partial#2}}
\def\deri#1#2{\frac{d #1}{d#2}}
\def\raiz{\sqrt{\tau}}
\def\flecha{\longrightarrow}
\def\fle{\rightarrow}
\def\gras{{\rm \ grad_{S^n}\ }}
\def\ds{\displaystyle}

\newtheorem{defi}{Definition}
\newtheorem{teor}{Theorem}
\newtheorem{prop}[teor]{Proposition}

\newtheorem{lema}[teor]{Lemma}
\newtheorem{lemap}{Lemma}

\newtheorem{nota}{Remark}

\newtheorem{corop}{Corollary}

\numberwithin{lemap}{teor}
\numberwithin{corop}{teor}
\numberwithin{ejer}{subsection}
\numberwithin{ejemplo}{subsection}
\begin{document}

\pagestyle{myheadings}\markboth{E. Cabezas-Rivas and V. Miquel}{Volume-preserving M.C.F. in Hyperbolic Space}%\markleft{}

%\markright{}

\title{Volume preserving mean curvature flow in the Hyperbolic space}% \footnote{Mathematics Subject
%Classification(2000) 53C42, 52C21}}

\author{ Esther Cabezas-Rivas and
Vicente Miquel}

\date{}

\maketitle

\begin{abstract}
We prove: \lq\lq If $M$ is a compact hypersurface of the hyperbolic space, convex by horospheres and evolving by the volume preserving mean curvature flow, then it flows for all time, convexity by horospheres is preserved and the flow converges, exponentially, to a geodesic sphere". In addition, we show that the same conclusions about long time existence and convergence hold if $M$ is not convex by horospheres but it is close enough to a geodesic sphere.\end{abstract}

\section{Introduction and Main Results }

Given an immersion $X:M\flecha \overline M$ of a compact $n$-dimensional manifold $M$ into a $(n+1)$-dimensional Riemannian manifold $\overline M$, the {\it mean  curvature flow} of $X$ is the solution of the partial differential equation 
\be\label{mcf}
\parcial{X_t}{t} = - H_t N_t, \text{ with the initial condition } X_0 = X,
\ee 
where $N_t$ is the outward unit normal vector of the immersion $X_t$ and  $H_t$ is the trace of the Weingarten map $L_{-N_t} = -L_{N_t}$ of $X_t$ associated to $-N_t$ (then, $H_t$ is $n$ times the usual mean curvature with the sign which makes positive the mean curvature of a round sphere in $\re^{n+1}$). From now on, by 
$M_t$ we shall denote both the immersion $X_t:M\flecha \overline M$ and the image $X_t(M)$, as well as the Riemannian manifold $(M,g_t)$
 with the metric $g_t$ induced by the immersion. The $n$-volume of $M_t$ (from now on called {\it area of }$M_t$) decreases along this flow, but no geometric invariant is preserved along it. 
 
 A related flow is the {\it volume preserving mean curvature flow}, which is defined as a solution of the equation
\begin{equation}\label{vpmf}
\parcial{X_t}{t} = (\oH_t- H_t)\ N_t,
\end{equation}
where $\oH_t$ is the averaged mean curvature
\be\label{oH}
\oH_t =\ds\frac{ \int_{M_t} H_t dv_{g_t} }{\int_{M_t} dv_{g_t} },
\ee
being $ dv_{g_t}$  the volume element on $M_t$. This flow also decreases the area of $M_t$, but preserves the volume of the domain $\Omega_t$ enclosed by $M_t$ (when such $\Omega_t$  exists).

In \cite{Hu84}, G. Huisken proved that, when $\overline M$ is the Euclidean space $\re^{n+1}$ and $M_0$ is strictly convex, then \eqref{mcf} has a maximal solution on a finite time interval $[0,T[$ and $M_t$ converges to a point as $t\to T$. Moreover, after appropriate rescaling of $X_t$ and $t$, $M_t$ converges to a round sphere. In 
\cite{Hu86}, he extended this result to compact hypersufaces in general Riemannian manifolds (with suitable bounds on curvature). 

The flow \eqref{vpmf} was considered by Huisken in \cite{Hu87}, again for $\overline M = \re^{n+1}$ and $M_0$ strictly convex; he proved that \eqref{vpmf} has a solution on $[0,\infty[$, which stays convex all time and  converges to a round sphere. However, he noticed the difficulty that the presence of averaged mean curvature in \eqref{vpmf} causes in order to extend this result to general Riemannian manifolds. In fact, Huisken illustrated this by showing a way to obtain examples of convex hypersurfaces in the sphere $S^{n+1}$ which could lose convexity along the flow. The idea is that,  if a piece $M'$ of $M_0$ is a part of a geodesic sphere of mean curvature near to $0$ and far from $\oH_0$, $M'$ moves in the outward radial direction of $M'$, soon becoming a totally geodesic hypersurface and, after that, changing the sign of the mean curvature. 

The above remark, pointed by Huisken in \cite{Hu87}, was really inspiring for us. First we noticed that examples like those in $S^{n+1}$ cannot happen in the Euclidean space, because when a geodesic sphere is moving outward in the the direction of its radius, it becomes of lower and lower normal curvature, but it never becomes a totally geodesic submanifold. Nevertheless, the Euclidean case was already settled by Huisken in \cite{Hu87}.

On the other hand, we realized that  a similar situation happens in the hyperbolic space: when a geodesic sphere moves outward in the radial direction, its normal curvature decreases, and it becomes nearer and nearer to that of a horosphere (see Remark (ii) below for a definition), but it never gets the curvature of a horosphere. The former intuitive idea was indeed the detonating clue which leads us to hope for a theorem like that of Huisken in \cite{Hu87} for the volume preserving mean curvature flow in the  hyperbolic space  of a   hypersurface convex by horospheres. This paper achieves the realization of such hope by proving the following theorem.

\begin{teor}\label{math}
Let $\eme$ be the complete simply connected $(n+1)$-dimensional hyperbolic space of sectional curvature $\lambda<0$. If $M_0$ is a compact hypersurface  convex by horospheres, then the equation  \eqref{vpmf} with initial condition $M_0$ has a unique solution $M_t$ such that
\begin{enumerate}
\item[(a)] it is defined for $t\in [0,\infty[$, 

\item[(b)] the hypersurfaces $M_t$ stay smooth and convex by horospheres for all time, 

\item[(c)] and the $M_t$'s converge exponentially (as $t \to \infty$, in the $C^m$ topology for any fixed $m\in\mathbb N$) to a geodesic sphere of $M_\lambda^{n+1}$ enclosing the same volume as $M_0$.
\end{enumerate}
\end{teor}

Next we include some remarks for a better understanding of the above statement.

{\bf Remarks} \ (i) Recall that a {\it horosphere} of $\eme$ is a hypersurface $\H$ obtained as the limit of a geodesic sphere of $\eme$ when its center goes to the infinity along a fixed geodesic ray, which is equivalent to say that $\H$ is a complete embedded hypersurface with normal curvature $\ral$. An {\it horoball} is the convex domain which boundary is an horosphere. 

(ii) A
hypersurface $M$ of $\eme$ is called {\it convex by horospheres} ({\it $h$-convex} for short) if  it bounds a domain
 $\Omega$ satisfying that, for every $p\in M=\partial\Omega$, there is a horosphere $\H$
of $\eme$ through
$p$ such that $\Omega$ is contained in the horoball of $\eme$ bounded by
$\H$. This $\H$ is called a {\it supporting horosphere} of $\Omega$ (and
of
$M$) through $p$. One shows that a hypersurface $M$ of $\eme$ is $h$-convex  if and only if all its normal curvatures are bounded from below by $\ral$. 

(iii) Usually two immersions $X_i:M_i \flecha \overline M$, $i=1,2$ are considered equivalent if there is a diffeomorphism $\phi:M_1\flecha M_2$ such that $X_1=X_2\circ \phi$. In this case, $X_1$ and $X_2$ are called parametrizations of the same immersed submanifold. For this reason, dealing  with submanifolds, one says that an immersion $X:M_2\flecha \overline M$ is in a neighborhood $U$ of $Y: M_1 \flecha \overline M$ in the $C^k$-topology if there is a diffeomorphism $\phi$  such that $X\circ\phi \in U$.  In  Theorem \ref{math}, we use the convergence in the $C^k$-topology in this sense.

(iv) Let us notice that $M_0$ $h$-convex implies that it is diffeomorphic to a sphere, and this implies that  $X$ is, in fact, an embedding. 

About the techniques applied for proving Theorem \ref{math}, as in  \cite{Hu87}, we use essentially maximum principles, employing also ideas  used by B. Andrews (\cite{An01}) and J. McCoy (\cite{MC04})  for similar theorems in the Minkowski (normed)  and Euclidean spaces, respectively. It may seem a surprising fact that the results in \cite{Hu87}, published in 1987, have not been exported previously to the hyperbolic space, because (as we remark above) they are based on a very natural idea from elementary hyperbolic geometry.  The reason for such delay could be that the proof of Theorem \ref{math} requires  an additional ingredient: a deep knowledge and a strong use of the geometric properties of the $h$-convexity described in the papers \cite{BGR}, \cite{BM99}, \cite{BM02} and \cite{BV}. 

Moreover, a method described in \cite{EsSi98} is used to prove that the convergence is exponential. This method relies on maximal regularity theory and is of independent interest. Indeed, its strength  allows us to extend statements (a) and (c) in Theorem \ref{math} to certain non-necessarily $h$-convex initial data. With more precision, as a by-product of the proof of the exponential convergence in Theorem \ref{math}, we shall obtain

\begin{teor}\label{corm} Let $\mc S$ be a geodesic sphere of $\eme$ and $0 < \beta < 1$. There exists an $\varepsilon>0$ such that, for every embedding  $X:M\flecha \eme$ with $h^{1 + \beta}$-distance to $\mc S$ lower than $\varepsilon$, the  equation \eqref{vpmf} has a unique solution satisfying $X_0=X$, defined on $[0,\infty[$ and which  converges exponentially to a geodesic sphere in $M_\lambda^{n+1}$ $h^{1 + \beta}$-close to $\mc{S}$ and enclosing the same volume as $X(M)$.
\end{teor}
$h^{1+\beta}({\mc M})$ denotes, for a compact manifold $\mc M$, the little Hšlder space of order $1+\beta$, that is, the closure of $C^\infty({\mc M})$ in the usual Hšlder norm of $C^{1+\beta}(M)$.

In a recent paper (cf.\cite{AlFr03}), Alikakos and Freire  proved long time existence for solutions of \eqref{vpmf} and convergence to constant mean curvature hypersurfaces in general ambient manifolds $\overline M$, but with the hypotheses that the initial condition  $M_0$ is \lq\lq close enough'' to a geodesic sphere of $\overline M$ (although $M_0$ does not need to be convex) and the scalar curvature of $\overline M$ has nondegenerate critical points. It may seem that such result includes our Theorem \ref{corm}, but this is not the case because $\eme$ has constant scalar curvature. 

The paper is organized as follows. In section \ref{notation}, we establish some notation and summarize the basic inequalities for $h$-convex sets which will be  used all along the proof of Theorem \ref{math}.  In section \ref{convexity}, we shall prove that the solution $M_t$ remains $h$-convex along all time it exists.  Section \ref{boundH} contains the main part of the proof:  the obtaining  of an universal (not depending on $t$) bound for $H_t$ and all its derivatives. As a consequence of this, we get that $M_t$ exists for $t\in [0,\infty[$. Sections 5 and 6 are devoted to prove statement (c) in Theorem \ref{math}: first, we find a time sequence  $\{t_i\}$ such that $\{M_{t_i}\}$ converges (up to isometries) to a geodesic sphere in $M^{n+1}_\lambda$; later, in section 6, we conclude that the full family $\{M_t\}$ converges $C^m$-uniformly and at exponential rate. Finally, the proof of Theorem \ref{corm} is included in section 7.

\section{Notation and preliminaries on $h$-convex sets}\label{notation}

From now on,  $\<\cdot,\cdot \>$, $\nablaa$, $\Deltaa$ and $\grada$ will denote
the metric, the covariant derivative, the Laplacian and the gradient
(respectively) of the ambient manifold $\eme$.  For
$\Deltaa$ (and the analog rough laplacian on tensor fields) we shall use the following  sign convention:
$$\Deltaa f =  \tr \nablaa^2 f.$$
The corresponding operators on $M$ will be denoted by
$\nabla$, $\Delta$ and $\grad$. 

When $\lambda<0$, we shall use the  notation: 

$$\s(t)= 
\frac{\sinh(\sqrt{|\lambda|}
t)}{\sqrt{|\lambda|}},  \quad
 \c(t)=\s'(t),\quad 
\ta(t) = \frac{\s(t)}{\c(t)}, \quad \text{and} \quad
\co(t) =
\frac{\c(t)}{\s(t)}.
$$ 
The functions above satisfy the following computational rules:
\begin{equation}\label{trirel}
c_\lambda^2 + \lambda\ s_\lambda^2=1, \quad c_{4\lambda} =
c_\lambda^2 - \lambda\ s_\lambda^2, \quad s_{4\lambda}=
s_\lambda c_\lambda.
\end{equation}

Given any point $p$ in the ambient space $\eme$, we
shall denote by $r_p$ the function \lq\lq distance to
$p$'' in $\eme$. Given a function $f: \re \flecha \re$, $f(r_p)$ will mean $f\circ r_p$. We shall also use the
notation $\partial_{r_p} = \grada r_p$. In the following lemma, we recall some formulae involving derivatives of $f(r_p)$ that we shall apply later. 

\begin{lema}[\cite{CGM}, \cite{GW},
\cite{Pe}]\label{lema5} In $\eme$, 
\begin{align}\label{compHR}
\displaystyle
&\<\nablaa_X\partial_{r_p},Y\> = \nablaa^2 r_p(X,Y) =  \begin{cases}  0 &{\rm if\ } X=
\partial_{r_p}\\
  co_\lambda(r_p)
\<X , Y\> &{\rm if\ } \<X,\partial_{r_p}\>=0
\end{cases} ,\\
&\quad \Deltaa r_p = n\ co_\lambda(r_p). \label{compHR2}
\end{align}
Moreover, if $f:\re\flecha\re$ is a $C^2$ function, 
\begin{align}\label{odeltafr}
&\overline\Delta (f(r_p)) =  f''(r_p) + f'(r_p)\  \overline\Delta r_p.
\end{align}
And, for the restriction of $r_p$ to a hypersurface $M$ of $\eme$, one has
\begin{align}\label{deltar} 
\Delta r_p &= -H \<N, \partial_{r_p}\> +  \co(r_p) \left( n\ -   |\partial_{r_p}^\top|^2\right). \\
 \Delta (f(r_p)) &= f''(r_p)\   |\partial_{r_p}^\top|^2 + f'(r_p)\   \Delta r_p \label{deltafr} \\
& \qquad = (  f''(r_p) - f'(r_p) \  \co(r_p)) |\partial_{r_p}^\top|^2  \nonumber\\
& \qquad \quad +  \ f'(r_p)\  ( n \ \co(r_p)    - H \<N, \partial_{r_p}\>). \nonumber
\end{align}
\end{lema}
Here $\partial_{r_p}^\top$ is the component of $\partial_{r_p}$ tangent to $M$, and it satisfies $\partial_{r_p}^\top = \grad(r_p\left|_M\right.)$

Next theorem summarizes some results contained in the quoted references.

\begin{teor}[\cite{BGR}, \cite{BM99}, \cite{BM02} and \cite{BV}]\label{TBGRM}  Let $\Omega$ be a compact $h$-convex domain and let $o$ be
the center of an inball of $\Omega$. If $\rho$ is the inradius of $\Omega$ and
$\tau=\ta\frac{\rho}{2}$, then 
\begin{itemize}
\item [a)] the maximal distance $\mad(o,\partial\Omega)$
between
$o$ and the points in $\partial\Omega$ satisfies the inequality
$$\mad(o,\partial\Omega)\leq \rho + \ral \ln\frac{(1+\raiz)^2}{1+\tau} < \rho + \ral \ln
2.$$
\item[b)] For any interior point $p$ of $\Omega$, $\<N,\partial_{r_p}\>  \ge \ral\  \ta(\dist(p,\partial\Omega))$, where $\dist$ denotes the distance in the ambient space $\eme$.
\end{itemize}
\end{teor}

\medskip
Moreover, in section \ref{convergence} we shall use the  elementary result stated below.

\begin{prop}\label{minproy}
In the Euclidean space $\re^{n+1}$, let $N\ne \zeta$ be two unit vectors. The maximal value of the (acute) angle between   a vector $v$ in the vector hyperplane $N^\bot$ orthogonal to $N$ and its projection onto the hyperplane $\zeta^\bot$ orthogonal to $\zeta$ is attained at the vectors in the intersection line of $N^\bot$ and the plane generated by $N$ and $\zeta$.
\end{prop}

We finish this section recalling the following consequence of the inequality between the trace and the norm of an endomorphism that will be used  through this paper:

\be\label{H2lenL2}
 |\nabla^m H|^2 \le  n \  |\nabla^m L|^2 \quad \text{ for every }  \quad m=0,1,2,3, ... 
\ee

%%%%%%Seccion 3
\section{Preserving $h$-convexity}\label{convexity}
\noindent With the notations of Theorem \ref{math},  here we shall prove
\begin{prop}\label{preserhc} In $\eme$, if  $M_0$ is  $h$-convex, under the volume preserving mean curvature flow \eqref{vpmf}, $M_t$ remains  $h$-convex for all the time such that the solution exists.
\end{prop}
For the proof of this result, we shall use the maximum principle for symmetric tensors as it is stated in
\cite{ChKn}, page 97.
Before, we need some evolution equations.

\begin{lemap}\label{l.eveq} For an arbitrary ambient space $\overline M$, the evolution equations of the metric $g_t$ and the second fundamental form $\a_t$  of a solution $M_t$ of \eqref{vpmf} are
\begin{align}
\parcial{g_t}{t} &= 2 (\overline H_t -  H_t)  \a_t  , \label{evol.g}\\
\parcial{\a_t}{t} &= \Delta_t \a_t  - 2 \ H_t  \<L_t^2\ \cdot\  ,\ \cdot\ \> + \oH_t  \left(\< L_t^2\ \cdot\  ,\ \cdot\ \> - \oR_{\ \cdot\ N_t\ \cdot\ N_t} \right)\label{evol.a}\\ 
 &+ (|L_t|^2 + \oRic(N_t,N_t))  \a_t  
 -  \oRic(\cdot\  ,\ L_t\  \cdot\ ) \nonumber\\
 &-  \oRic(L_t \ \cdot\  ,\cdot\ ) +  \oR(N_t,  L_t \  \cdot\  , N_t,\ \cdot ) \nonumber \\
&  + \oR(N_t,  \  \cdot\  , N_t, L_t \  \cdot ) + 2 RiL_t - \overline\nabla_\cdot \oRic(\cdot, N_t) - \overline\delta\oR_{\ \cdot\  \cdot\  N_t},  \nonumber
 \end{align}
 where $\oR$ and $\oRic$ denote, respectively, the curvature and  Ricci tensors of $\overline M$, $RiL_t(Z,Y)= \sum_{i=1}^n  \oR_{e_i,Z,L_te_i, Y}$  and $ \overline\delta\oR= 
\sum_{i=1}^n  \overline \nabla_{e_i}(\oR)_{e_i}$ for some local orthonormal frame $\{e_i\}$ of $M_t$.
\end{lemap}
\begin{demo} Formula \eqref{evol.g} follows from \eqref{vpmf} by a direct computation as in \cite{Hu84}. Also by this way one obtains 
\begin{align}\label{evol.a0}
&\parcial{\a_t}{t} = -(\oH_t-H_t) \oR_{\ \cdot\ N_t\ \cdot N_t} + \nabla^2H_t + (\oH_t - H_t) \<L_t^2\ \cdot\  ,\ \cdot \>,
\end{align}
and, having into account the (generalized)   Simons' formula for the rough Laplacian of the second fundamental form (see, for instance, \cite{Ber})
\begin{align}
\Delta \alpha &= \nabla^2 H + H \<L^2 \ \cdot\ , \ \cdot\ \> - |L|^2 \a + H \oR_{N\ \cdot\ N\ \cdot\ } - \oRic(N,N)\  \a \label{Simons}\\
&  - 2 \  RiL +   \oRic(\cdot\  ,\ L \ \cdot\ ) +  \oRic(  L \  \cdot\  ,\ \cdot\ ) \nonumber \\
& -  \oR(N,  L\  \cdot\  ,N, \ \cdot\ ) - \oR(N,  \  \cdot\  ,N, L \ \cdot\ ) + \overline\nabla_\cdot \oRic(\cdot, N) + \overline\delta\oR_{\ \cdot \  \cdot \  N}, \nonumber \end{align}
we get \eqref{evol.a}.
\end{demo}

\noindent
\begin{demo} {\bf of Proposition \ref{preserhc}}. Let us take $A_t = \alpha_t - \ral \ g_t$. Notice that, from the explicit expression of the  curvature tensor $\oR$ of $ \eme$, the equation \eqref{evol.a} becomes
\begin{align}
&\parcial{\a_t}{t} = \Delta_t \a_t + (\oH_t - 2 H_t)\  \a L_t + (|L_t|^2 - \lambda n)\  \a_t + \lambda (2 H_t - \overline H_t)\  g_t,\label{evol.am} 
\end{align}
where $ \a L$ is defined by $\a L(X,Y) = \a(L X, Y) $.

\noindent From \eqref{evol.g} and \eqref{evol.am}, we obtain
\begin{align}
&\parcial{A_t}{t} = \Delta_t A_t + B_t \text{, with } \label{evol.al}\\
&B_t = (\overline H_t - 2 H_t) (\a L _t-  \lambda  g_t) + \left(|L_t|^2 - \lambda n - 2 \ral (\overline H_t - H_t)\right)  \a_t. \nonumber
\end{align}
Let $V$ be a unitary null vector of $A_t$, that is, $L_t V = \ral V$. A straightforward computation gives
\begin{align*}
B_t(V,V) &=  \ral |L_t|^2 + 2 \lambda H_t - n \lambda \ral \\
& \ge  \frac{\ral}{n} H_t^2 + 2 \lambda H_t - n \lambda \ral = \frac{\ral}{n} (H_t - n  \ral )^2 \ge 0 ,
\end{align*}
using \eqref{H2lenL2} in the first inequality. Now, the proposition follows by the maximum priciple for symmetric tensors quoted above.
\end{demo}

\section{Long Time Existence}\label{boundH}

Along this section, we shall denote by $[0,T[$ the maximal interval where the solution of \eqref{vpmf} is well defined, and want to prove that $T=\infty$.

 The main point to establish long time existence is to show that  $|L_t|$ has an uniform bound independent of $t$. Since we proved previously that $M_t$ is $h$-convex as long as it exists, then $|L_t|^2 \le H_t^2$; therefore, it is enough to show that $H_t$ has an upper bound independent of $t$.  In order to achieve this, first we shall study the evolution under \eqref{vpmf} of the function
\be\label{defW}
W_t = \frac{H_t}{\sigma_t - c}, \text{ being } \sigma_t = s_\lambda(r_p)  \<N_t, \partial_{r_p}\>  \text{ and }  c \text{ any constant. }\ee

 Before starting the way to obtain the evolution equation for $W_t$, we would like to remark that $\sigma_t$ depends on the choice of the point $p$. This fact will be important later, when we write inequalities. 

Let $g_t^\flat$ denote the metric induced on $T^*M$ by $g_t$ through the isomorphism $\flat_t:TM\flecha T^*M$ of lowering indices. The matrix $(g^{ij})$ of $g_t^\flat$ in some basis is the inverse of the matrix $(g_{ij})$ of $g_t$ in the dual basis. Using this fact and \eqref{evol.g}, one obtains
\be\label{evol.gf}
\parcial{g_t^\flat}{t} = - 2 (\overline H_t -  H_t)  \a^\flat_t , 
\ee
where $\a_t^\flat$ denotes, again, the tensor induced on $T^*M$ by $\a_t$ through the isomorphism $\flat_t$. From \eqref{evol.a0} and \eqref{evol.gf}, we get
\be\label{evol.H}
\parcial{H_t}{t} = \Delta_t H_t + (|L_t|^2 + \oRic(N_t,N_t) ) (H_t-\oH_t),
\ee
which, by the expression of the curvature tensor $\oR$ of $\eme$, becomes
\be\label{evol.He}
\parcial{H_t}{t} = \Delta_t H_t + (|L_t|^2 + n\ \lambda ) (H_t-\oH_t) .
\ee
Another standard computation (similar to that done in \cite{Hu84}) allows to obtain, from \eqref{vpmf}, the evolution equation
\be\label{evol.N}
\frac{\nablaa N_t}{\partial t} = - \grad(\oH_t - H_t) = \grad H_t.
\ee

Now, let us note that, for any smooth function $\varphi:\re\flecha \re$, a direct calculation using \eqref{compHR} and \eqref{vpmf} gives 
\be\label{evol.phr}
\Daparcial{(\varphi(r_p) \partial_{r_p})}{t} =(\oH_t-H_t) \left((\varphi'(r_p) - \varphi(r_p) \ \co(r_p)) \<\partial_{r_p}, N_t\> \partial_{r_p} + \varphi(r_p )\  \co(r_p)\  N_t\right).
\ee
Taking $\varphi = \s$ and using \eqref{evol.N}, we arrive to
\be\label{evol.sig}
\parcial{\sigma_t}{t} = (\oH_t-H_t)\  \c(r_p) + \s(r_p)  \< \partial_{r_p}, \grad H_t\>.
\ee
From  \eqref{deltafr} (with $f=\s$) and \eqref{deltar}, we have
\be\label{deltasr}
\Delta_t \left(\s(r_p)\right) = -\frac1{\s(r_p)} |\partial_{r_p}^\top|^2 - \c(r_p)\  H_t \<N_t, \partial_{r_p}\> + n \frac{\c^2}{\s}(r_p).
\ee
Straightforward computations having into account \eqref{compHR} give
\begin{align}
&\<\grad \ \s(r_p), \grad \<\partial_{r_p}, N_t \>\> \label{gragra} \\
& \qquad = -\frac{\c^2}{\s}(r_p) \<\partial_{r_p}, N_t\> |\partial_{r_p}^\top|^2 + \c(r_p) \ \a(\partial_{r_p}^\top, \partial_{r_p}^\top), \nonumber
\\
&\Delta_t \< \partial_{r_p}, N_t \> =   \frac1{\s^2(r_p)}\<\partial_{r_p},N_t\> |\partial_{r_p}^\top|^2\  - n\ \co^2(r_p) \<\partial_{r_p},N_t\>   \label{delta<>} \\
&\qquad \quad + \co(r_p) \<N_t,\partial_{r_p}\>^2  H_t  - 2\  \co(r_p)\ \a(\partial_{r_p}^\top,\partial_{r_p}^\top) \nonumber \\
& \qquad \qquad+ 2\   \co^2(r_p)\  \<N_t,\partial_{r_p}\> |\partial_{r_p}^\top|^2 \nonumber \\
& \qquad \quad + \<\partial_{r_p}^\top, \grad H_t\> + \co(r_p) \  H_t - \<\partial_{r_p}, N_t\> |L_t|^2. \nonumber
\end{align}
Joining \eqref{deltasr}, \eqref{gragra} and \eqref{delta<>}, we reach
\be\label{deltasigma}
\Delta_t \, \sigma_t = \c(r_p) \  H_t + \< \s(r_p)\ \partial_{r_p}, \grad H_t\> - \sigma_t \, |L_t|^2.
\ee
By substitution of this expression in \eqref{evol.sig}, we obtain the evolution equation 
\be\label{evol.sigD}
\parcial{\sigma_t}{t} = \Delta_t \sigma_t +  |L_t|^2 \sigma_t  + (\oH_t- 2\ H_t)\  \c(r_p) .
\ee

From \eqref{defW}, \eqref{evol.He} and \eqref{evol.sigD}, it follows
\begin{align}
\parcial{W_t}{t} &=  \frac{1}{\sigma_t - c}\  \Delta_t H_t + \frac{1}{\sigma_t - c}\  (H_t - \oH_t)\  (|L_t|^2 + n \lambda) \label{evol.W}\\
&-  \frac{H_t}{(\sigma_t - c)^2} \left( \Delta_t \sigma_t + |L_t|^2 \sigma + (\oH_t - 2\ H_t)\ \c(r_p)\right).
 \nonumber
 \end{align}
Taking definition \eqref{defW} as starting point,  another computation leads to 
\begin{align}
\Delta_tW_t &=  \frac{1}{\sigma_t - c}\  \Delta_t H_t +  \frac{2 H_t}{(\sigma_t - c)^3} \   |\grad \ \sigma_t|^2 - \frac{H_t}{(\sigma_t - c)^2}\   \Delta_t\sigma_t \label{DeltaW}\\
& - 2 \<\grad H_t,     \frac{1}{(\sigma_t - c)^2} \, \grad \ \sigma_t\>. \nonumber 
\end{align}
Replacing \eqref{DeltaW} into \eqref{evol.W} and doing a few more computations, we can write
\begin{align}
\parcial{W_t}{t} &= \Delta_t W_t + \frac{2}{\sigma_t - c} \<\grad W_t, \grad \, \sigma_t\> \label{evol.WD}\\
&-  \frac{\oH_t}{\sigma_t-c}\  (|L_t|^2 + n \lambda )  
- \frac{W_t}{\sigma_t-c} \   \oH_t\   \c(r_p)  + 2\  W_t^2\  \c(r_p) \nonumber \\
& - \frac{c}{\sigma_t - c}\  W_t \  |L_t|^2  + n\  \lambda\  W_t. \nonumber
\end{align}

To get fine and independent of $t$ bounds for $W_t$ from \eqref{evol.WD} by application of the maximum principle, previously we need to bound $r_p$ and $\<N_t,\partial_{r_p}\>$. In order to do so, we shall use Theorem \ref{TBGRM}.

\begin{lema} Let $\psi$ be the inverse of the function $\ds s\mapsto \vle(S^n) \int_0^s \s(\ell) d\ell$ and $\xi$ the inverse function of $s\mapsto s+ \ds  \ral \ln\frac{(1+\sqrt{\ta(\frac{s}{2})})^2}{1+\ta(\frac{s}{2})} $. If $V_0=\vle(\Omega_0)$ and $\rho_t$ is the inradius of $\Omega_t$, then
\be\label{brho}
\xi(\psi(V_0)) \le \rho_t \le \psi(V_0),
\ee
for every $t\in [0,T[$.
\end{lema}
\begin{demo}
Since the flow preserves the enclosed volume, we have $\vle(\Omega_t) = V_0$ for all $t\in [0,T[$. If we take spherical geodesic coordinates in $\eme$ around a center $p_t$ of an inball of $\Omega_t$, we can describe $M_t$ as the graph of a function $\ell:S^n \flecha \re^+$, and the volume of $\Omega_t$ is given by
\be
\vle(\Omega_t) = \int_{S^n} \int_0^{\ell(u)} \s^n(s)\  ds \  du.
\ee
But $\rho_t \le \ell(u) \le \mad(p_t,M_t) \le \rho_t + \ds  \ral \ln\frac{(1+\sqrt{\ta(\frac{\rho_t}{2})})^2}{1+\ta(\frac{\rho_t}{2})} $ (where we have used Theorem \ref{TBGRM} a) for the last inequality), thus the lemma follows having into account that $\psi^{-1}$ and $\xi^{-1}$ are increasing functions.
\end{demo}

An immediate consequence of the lemma above and Theorem  \ref{TBGRM} a) is

\begin{corop}\label{dpqbo}
For every $t\in [0,T[$, if $p,q \in \Omega_t$, then 
\be\label{dpqboe} \dist(p,q) < 2 (\psi(V_0) +  \ral \ \ln 2).\ee
\end{corop}

Now, let us continue with the task of bounding $W_t$. First, we fix an arbitrary $t_0 \in [0,T[$.  As before, $\rho_t$ will denote  the inradius of  an inball of $\Omega_t$ and $p_t$ its center.  Although $p_{t_0}$ does not need to be the center of an inball of $\Omega_t$ when $t\ne t_0$, we can use Corollary \ref{dpqbo} whereas $p_{t_0}\in\Omega_t$ to bound $r_{p_{t_0}}(x) \le 2 (\psi(V_0)+  \ral \ \ln 2)$ for every $x\in M_t$. Consequently, our next goal is to estimate a time interval $[t_0,t_0+\tau[$ such that $p_{t_0}\in\Omega_t$ for $t\in [t_0,t_0+\tau[$. To do so,  we shall compare the motion of $M_t$ following the equation \eqref{vpmf} with the motion  under \eqref{mcf} of a geodesic sphere centered at $p_{t_0}$ with radius $\rho_{t_0}$ at time $t_0$. As a result, we shall obtain

\begin{lema}\label{t0+epsi}
There is $\tau=\tau(\lambda, n, V_0) >0$ such that $B(p_{t_0},\rho_{t_0}/2) \subset \Omega_t$ for every $t \in [t_0, t_0+\min\{\tau, T-t_0\}[$.
\end{lema}
\begin{demo}
Let $r_B(t)$ be the radius  at time $t$ of a geodesic sphere $\partial B(p_{t_0},r_B(t))$ centered at  $p_{t_0}$, evolving under \eqref{mcf} and with the initial condition $r_B(t_0) = \rho_{t_0}$. From \eqref{mcf}, \eqref{compHR} and the fact that the mean curvature of a geodesic sphere centered at $p_{t_0}$ is $\Deltaa r_{p_{t_0}}$, we get
\be\label{evol.rB}
\ds \parcial{r_B(t)}{t} = - n\   \co (r_B(t)),
\ee
 and the solution of this differential equation satisfying $r_B(t_0)= \rho_{t_0}$ is
\be\label{rBt}
\c(r_B(t)) = e^{\lambda n (t-t_0)} \c(\rho_{t_0}).
\ee
Then, for $t\ge t_0$ (and because $\c$ is an increasing function), $r_B(t) \ge \rho_{t_0}/2$ if and only if $e^{\lambda n (t-t_0)} \c(\rho_{t_0}) \ge \c(\rho_{t_0}/2)$, i.e., 
\be
\text{$r_B(t) \ge \rho_{t_0}/2$ \quad  if and only if \quad $t-t_0 \le \ds \frac1{- \lambda\ n} \ln \frac{ \c(\rho_{t_0}) } {\c(\rho_{t_0}/2)}$  } \nonumber
\ee 
and, as the function $\ds s \mapsto \ln \frac{ \c(s) } {\c(s/2)}$ is increasing, using \eqref{brho}, we  have
\be\label{rBr2}
\text{$r_B(t) \ge \rho_{t_0}/2$ \quad if \quad $t-t_0 \le \ds  \frac1{- \lambda\ n} \ln \frac{ \c(\xi(\psi(V_0)))}{\c(\xi(\psi(V_0))/2)}=:\tau$}.
\ee 

For any $x\in M$, let $r(x,t)= r_{p_{t_0}}(X_t(x))$. From \eqref{vpmf}, it follows 
\be \label{evol.rx}
\parcial{r}{t} = (\oH_t - H_t) \< N_t, \partial_{r_{p_{t_0}}}\>.
\ee
If $\varphi:\re\flecha\re$ is a function satisfying $\varphi'(s)=\ta(s)$, and we set $f(x,t) = \varphi(r(x,t)) -\varphi(r_B(t))$, from \eqref{evol.rB} and \eqref{evol.rx}, we obtain
\be\label{f't}
\parcial{f}{t} = \ta(r_{p_{t_0}}) \  (\oH_t - H_t) \<N_t,\partial_{r_{p_{t_0}}}\> + n.
\ee
On the other hand, from \eqref{deltafr},
\begin{align}\label{Deltaphir}
\Delta_t f & = \Delta_t (\varphi(r_{p_{t_0}})) \\
&=  \left( \frac1{\c^2}(r_{p_{t_0}}) - 1\right) |\partial_{r_{p_{t_0}}}^\top|^2 +   n  - H_t\ \ta(r_{p_{t_0}}) \<N_t, \partial_{r_{p_{t_0}}}\> .\nonumber
\end{align}
Now, let $t_1= \inf\{t>t_0;\ p_{t_0}\notin\Omega_t\}$. Because $\Omega_t$ is $h$-convex, $\<N_t, \partial_{r_{p_{t_0}}} \> \ge 0 $ for $t\in[t_0,t_1]$. By substitution of \eqref{Deltaphir} into \eqref{f't}, we arrive to
\begin{align}
\parcial{f}{t} = \Delta_t f & +\ta(r_{p_{t_0}}) \  \oH_t  \<N_t,\partial_{r_{p_{t_0}}}\>  \label{f'tDci} \\
&+ \left( 1 - \frac1{\c^2}(r_{p_{t_0}}) \right) |\partial_{r_{p_{t_0}}}^\top|^2  \ge \Delta_t f, \label{f'tD} \quad \text{ and } \nonumber\\
    f(x,t_0) = &\varphi(r(x,t_0)) -\varphi(\rho_{t_0})   \ge 0. 
\end{align}
Using the scalar maximum principle for parabolic inequalities (cf. \cite{ChKn}, page 94) gives $f(x,t)\ge 0$ for $t_0\le t\le t_1$ as long as $f(x,t)$ is well defined. But $r(x,t)$ is well defined for $0\leq t <T$, and it follows from \eqref{rBt} that $r_B(t)$ is well defined (that is,  positive) for $t\in [t_0, t_0 - \frac1{\lambda n} \ln(\c(\rho_{t_0}))[ \supset [t_0, t_0+\tau[$. Then $f(x,t) \ge 0$ on $[t_0,  \min\{t_0+\tau, T, t_1\}[$. 

Notice that, by definition of $t_1$, $p_0\in M_{t_1} =\partial \Omega_{t_1}$. If $t_1< \min\{t_0+\tau, T\}$, $f(x,t_1)\ge 0$ and $B(p_{t_0},r_B(t_1))\subset \Omega_{t_1}$, which is a contradiction; therefore, $t_1 \ge \min\{t_0+\tau, T\}$, and the lemma follows.
\end{demo}

As a consequence of \eqref{brho}, Corollary \ref{dpqbo} and Lemma \ref{t0+epsi}, on the interval $[t_0, t_0+\min\{\tau, T-t_0\}[$, and on the hypersurface $M_t$
$$C:=\frac{\xi(\psi(V_0))}{2} \le r_{p_{t_0}} \le 2 (\psi(V_0)+  \ral \ \ln 2)=:D.$$
Moreover, having into account  Theorem \ref{TBGRM} b),
$$\sigma_t = \s(r_{p_{t_0}}) \<N_t, \partial_{r_{p_{t_0}}}\> \ge \ral \  \s(C) \   \ta(C).$$
Then, if we take the constant $c$ in the definition \eqref{defW} as 
$$c=  \frac{\ral\ \s(C)\ \ta(C)}{2},$$ 
 we get $\sigma_t-c \ge c >0$.

Let us go back to equation \eqref{evol.WD}. From the above remark on $\sigma_t-c$ and the $h$-convexity of $M_t$,  we have $W_t\ge 0$ and $H_t + n\lambda \ge 0$. Moreover, $|L_t|^2 \ge \frac1n H_t^2$. Now we can use these inequalities in \eqref{evol.WD} to obtain
 \begin{align}
\label{evol.Wineq} \parcial{W_t}{t}  &\le \Delta_t W_t + \frac{2}{\sigma_t - c} \<\grad W_t, \grad \, \sigma_t\> + 2 \  \c(D) \ W_t^2 - \frac{c}{\sigma_t - c} W_t \frac{H_t^2 }{n} \\ &\le \Delta_t W_t + \frac{2}{\sigma_t - c} \<\grad W_t, \grad \, \sigma_t\> + 2\   \c(D) \ W_t^2 - \frac{c^2}{n} W_t ^3. \nonumber
\end{align}

By other version of the scalar maximum principle (cf. \cite{ChKn}, page 96),  in the interval  $[t_0, t_0+\min\{\tau, T-t_0\}[$,  $W_t(x)$ is bounded from above by the solution $w(t)$ of  the ordinary differential equation 
$$ w' = \left(2  \c(D)  - \frac{c^2}{n} w \right) w^2 , \quad \text{with} \quad w(t_0)= \max_{x\in M} W_{t_0}(x).$$
Observing that  $w'(t)<0$ when $w> \frac{2 n \c(D)}{c^2}$, it is straightforward  to show that  $w(t) \le \max\{w(t_0), \frac{2 n \c(D)}{c^2}\}$. Thus we deduce
$$W_t(x) \le  \max\{\max_{x\in M} W_{t_0}(x), \frac{2 n \c(D)}{c^2}\} \text{ for every $t\in [t_0, t_0+\min\{\tau, T-t_0\}[$}. $$

 From the definition of $W_t$, the election of  $c$ and the upper bound of $\rho_t$, we have 
$$H_t(x) \le (\s(D) - c) \max\{\max_{x\in M} W_{t_0}(x), \frac{2 n \c(D)}{c^2}\}.$$
Since this occurs for any $t_0$ and $\tau$ does not depend on $t_0$, we arrive to
\begin{align}\label{bound.H}
 & H_t(x) \le   (\s(D) - c) \max\{\max_{x\in M} W_{0}(x), \frac{2 n \c(D)}{c^2}\} =: C(n,\lambda,M_0) \\
 & \text{ for every } t\in [0,T[.  \nonumber
\end{align}
 This implies, by the definition of $\oH_t$ and the $h$-convexity of $M_t$,
\begin{align}\label{Lbound}
&\oH_t \le C(n,\lambda,M_0) \text { and } |L_t|^2 \le H_t^2 \le C(n,\lambda,M_0)^2 \\
&\text{ for every } t \in [0,T[. \nonumber
\end{align}

From \eqref{evol.a}, reasoning like in \cite{Hu84} and \cite{Ham82}\S 13, one can deduce, for every natural number $m$, the following evolution equation 
\begin{align*}
\parcial{}{t} |\nabla^m L_t |^2 = & \Delta_t |\nabla^m L_t|^2 - 2 |\nabla^{m +1} L_t |^2 + C(m, n, 
\lambda) |\nabla^m L_t |^2 
\\ &+ \sum_{i + j + k = m} \nabla^i L_t *\nabla^j L_t  * \nabla^k L_t * 
\nabla^m L_t 
 + {\oH_t} \sum_{i + j = m} \nabla^i L_t *\nabla^j L_t * \nabla^m L_t
\end{align*}
Then, using \eqref{Lbound} and arguing in the same way as in \cite{Hu87} Theorem 4.1, we conclude

\begin{prop} For every natural number $m$, there is a constant $C_m(n,\lambda,M_0)$ such that
\be |\nabla^mL_t|^2 \le C_m(n,\lambda,M_0) \label{Lmbound}  \ee
\end{prop}

From \eqref{evol.g}, \eqref{bound.H}, \eqref{Lbound}, and \eqref{Lmbound}, it follows (like in \cite{Hu84} pages 257, ff.) that, if $T<\infty$, then $X_t$ converges (as $t\to T$, in the $C^\infty$-topology) to a unique smooth  limit $X_T$ which represents a smooth $h$-convex hypersurface. Now we can apply the the short time existence theorem to continue he solution after $T$, arriving to a contradiction.
In short, the proof that the  solution of \eqref{vpmf} is defined on $[0,\infty[$ (that is, the long time existence statement in Theorem \ref{math}) is finished.

 \section{Convergence to a geodesic sphere}\label{convergence} 
 
Observe that, to finish the proof of Theorem \ref{math}, it remains to deal with the issues related to the convergence of the flow. We begin this task in the present section by proving 
 \begin{prop}\label{SeqConv}
 There is a sequence of times $t_{1} < t_2 < ... < t_{k} < ...\to \infty $ and isometries $\varphi_{t_1},\ \varphi_{ t_2},\ ... \ , \varphi_{ t_k},\ ...  \ $ of  $\eme$ such that $\varphi_{t_i}(M_{t_i})$ $C^\infty$-converges to an embedded geodesic sphere.
 \end{prop}
 \begin{demo} 
 The proof  is organized in two main steps. Let us begin by showing that if the aforementioned limit exists, it should be a hypersurface in $M^{n +1}_\lambda$ of constant mean curvature. As $H_t$ is invariant by the family of isometries $\{\varphi_t\}$, it will be enough to prove that $H_t$ (instead of $H_t \circ \varphi_t$) tends to a constant as $t\to \infty$, in other words,
 
 \medskip
 {\bf Step 1.} \label{LHtooH}
{\it The mean curvature $H_t$ of the hypersurfaces $M_t$ which evolve following  \eqref{vpmf} converges to its average, that is,
\be\label{HtooH}
\lim_{t\to \infty} \sup_{M_t} |H_t - \oH_t| = 0.
\ee}

 In order to prove the above claim, we shall state a series of auxiliary results.
 
 \medskip
 
$\bullet$ {\bf (\cite{Aub}, p.\,91)} {\it Let $M$ be a Riemannian manifold.  If a real function $f$ on $M$ satisfies $f\in L^1(M)$, $\int_M f dV = 0$ and $|\grad f| \in L^r(M)$, then
 \be \label{Aub1}
\sup_{M} | f| \le C\  ||\grad f||_r  \quad \text{ for every } r>n.
 \ee}

$\bullet$ {\bf (\cite{Aub}, p.\,93)} {\it Let $M$ be a Riemannian manifold and let $p$, $q$, $r$ be real numbers satisfying $1\le q, r \le \infty$ and $\ds\frac2{p} = \frac1{q} + \frac1{r}$. Every function $f\in C_0^\infty(M)$ satisfies
 \be\label{Aub2}
 ||\grad f||_p^2 \le (n^{1/2} + |p-2|)\ ||f||_q \   ||\nabla^2f||_r.
 \ee}
 
$\bullet$ {\bf (\cite{Aub}, p.\,89)} {\it Let $M$ be a Riemannian manifold and let $p$, $q$, $r$, $a$ be real numbers satisfying $1\le r < q  \le \infty$, $p\in[r,q]$ and $\ds a=\frac{1/p - 1/q}{ 1/r - 1/q}$. If $f\in L^r(M)\cap L^q(M)$, then $f\in L^p(M)$ and
 \be\label{Aub3}
 || f||_p \le  ||f||_r^a  \   ||f||_q^{1-a}.
 \ee}

Now we are in position to start proving \eqref{HtooH}.

From \eqref{evol.g}  and the expression of  $dv_{g_t}= \sqrt{\det(g_{t_{ij}})} \ du^1 ... du^n$ in local coordinates, a straigthforward computation gives
\be\label{evol.dV}
\parcial{}{t} dv_{g_t} =  (\oH_t - H_t) \ H_t \ dv_{g_t}.
\ee
This leads to
\be\label{evol.vol}
\deri{}{t} \vle(M_t) = - \int_M (H_t - \oH_t)^2 \ dv_{g_t},
\ee
thus 
\be\label{int0inf}
\int_0^\infty  \int_M (H_t - \oH_t)^2 \ dv_{g_t}\  dt = \lim_{t\to\infty} (\vle(M_0)-\vle(M_t)) \le \vle(M_0).
\ee

On the other hand, from \eqref{evol.dV}, \eqref{evol.He},   \eqref{Lmbound}  and \eqref{H2lenL2} it  follows that 
\linebreak 
$\ds \deri{}{t}\int_M(H_t-\oH_t)^2 dv_{g_t}$ is uniformly bounded. Therefore, \eqref{int0inf} implies
\be\label{limint}
\lim_{t\to \infty} \int_M (H_t-\oH_t)^2 dv_{g_t} = 0.
\ee

Since $\int_M (H_t-\oH_t) dv_{g_t} = 0$, $H_t-\oH_t$ is smooth and $M$ is compact,  then $H_t-\oH_t$ satisfies the hypotheses required to apply \eqref{Aub1}; so 
\be
\sup_{M_t}|H_t- \oH_t| \le C ||\grad(H_t-\oH_t)||_p \quad \text{ for every } p>n.
\ee

Using  \eqref{Aub3}, with $q=\infty$ and $r=2$, we get
\be \sup_{M_t}|H_t-\oH_t| \le C ||\grad(H_t-\oH_t)||_2^{2/p}  ||\grad H_t||_\infty^{1-2/p}. \label{59}
\ee

As a consequence of \eqref{H2lenL2} and \eqref{Lmbound}, one has the inequality
\be
|\grad H_t| \le   \sqrt{n \ C_1}. \label{510}
\ee

Moreover, if we apply \eqref{Aub2} to $f=(H_t-\oH_t)$, with $p=q=r=2$, we have
\be
||\grad (H_t-\oH_t)||_2^2 \le n^{1/2} ||H_t-\oH_t||_2 ||\nabla^2 H_t||_2 \label{511}
\ee

Replacing  \eqref{510} and \eqref{511} in \eqref{59}, we obtain that there is a constant $K$ depending only on $n$, $\lambda$ and $M_0$ such that 
\be \label{bHmh}
\sup_{M_t} |H_t- \oH_t| \le K\ (||H_t-\oH_t||_2\ ||\nabla^2H_t||_2)^{1/p}.
\ee
But, using again \eqref{H2lenL2} and \eqref{Lmbound}, and the decrease of $\vle(M_t)$ given by \eqref{evol.vol},
\begin{align}  \label{bD2H}
||\nabla^2 H_t||_2 &= \left(\int_M |\nabla^2 H_t|^2 \ dv_{g_t}\right)^{1/2} \le \sup_{M_t}|\nabla^2H_t| \  \vle(M_t)^{1/2}\\
&\le \sqrt{n} \  \sup_{M_t}  |\nabla^2L| \  \vle(M_0)^{1/2} \le K_1(n,\lambda,M_0).  \nonumber
\end{align}

By \eqref{limint}, \eqref{bHmh} and  \eqref{bD2H},  we reach
\begin{equation*}
\sup_{M_t} |H_t- \oH_t| \le K_2(n,\lambda,M_0)\ \left(\int_M (H_t-\oH_t)^2\ dv_{g_t}\right)^{1/(2p)}
\xrightarrow[t\to\infty]{} 0,
\end{equation*}
which finishes the proof of \eqref{HtooH}.

\medskip

Next step is to show the existence of the convergent sequence claimed in Proposition \ref{SeqConv}. With more precision, 
\medskip

{\bf Step 2.} {\it There exists a family of isometries $\{\varphi_t: M_\lambda^{n+1} \fle M_\lambda^{n+1}\}$ such that, if we consider the compositions $\varphi_t \circ X_t$ with $X_t$ being a solution of \eqref{vpmf} on $[0, \infty[$, then $\{\varphi_t \circ X_t : M \fle M_\lambda^{n+1}\}$ is precompact  in the $C^\infty$-topology. Moreover, the limit $\tilde M_\infty$ is a compact embedded hypersurface of $M_\lambda^{n+1}$}.

\medskip

 For each $t$, let us  fix a center $p_t$ of an inball of $\Omega_t$, and  let $\varphi_t$ be an isometry of $\eme$ carrying $p_t$ onto $p_0$. Obviously, each $\varphi_t(X_t(M))$ is an $h$-convex hypersurface with a center of an inball at $p_0$ and inradius $\rho_t$. Then, by Theorem \ref{TBGRM} and \eqref{brho}, $\dist(p_0, \varphi_t\circ X_t(x)) $ has an upper bound independent of $t$ and of $x$, i.e., the family $\{\varphi_t\circ X_t\}_{t\geq 0}$ is uniformly  bounded.
 
Let us denote by $S^n$ the unit sphere in $T_{p_0}\eme$.  For each $t$, since  $\varphi_t(X_t(M))$ is $h$-convex, there exists a function $\r:S^n \flecha \re^+$ such that we can parametrice $\varphi_t(M_t)$ by a map $\tilde X_t : S^n \flecha \eme$ satisfying 
\be\label{Xtilde}
\tilde X_t(u) = \exp_{p_0} \r(u) u.
\ee
Notice that $\r(u)= r_{p_0}(\tilde X_t(u))$. For any local orthonormal frame $\{e_i\}_{i=1}^n$  of  $S^n$, we have
\begin{align} \label{Xtei}
\tilde X_{t* \, u} e_i &= \exp_{p_0*} (e_i\r)(u) \ u +  \exp_{p_0*} \r(u)\ e_i \\
 &= e_i(\r)(u) \partial_{r_{p_0}} + \s(\r(u)) \tau_s e_i, \nonumber
\end{align}
where $\tau_s$ denotes the parallel transport along the geodesic starting from $p_0$ in the direction of $u$, and until $\exp_{p_0} \r(u) u$. 

Let $N_t$ be the outward unit normal vector to $\varphi_t(M_t)$. Observe that, by \eqref{Xtei}, the projection $\pi_\bot \tilde X_{t*} e_i$ of  $\tilde X_{t*} e_i$ onto the space $\partial_{r_{p_0}}^\bot$ orthogonal to $\partial_{r_{p_0}}$ is $\s(\r) \tau_s e_i$. Using Proposition \ref{minproy} (with $\zeta = \partial_{r_{p_0}}$ and $N = N_t$), the angle $\beta$ between $\tilde X_{t*}e_i$ and its projection is bounded from above by the angle $\beta_0$ they  form in case $\tilde X_{t*}e_i$, $\partial_{r_{p_0}}$ and $N_t$ are in the same plane. Then, in general, 
$$\s(\r) =  |\pi_\bot \tilde X_{t*} e_i| =  |\tilde X_{t*} e_i| \cos \beta \ge   |\tilde X_{t*} e_i|  \cos\beta_0 = |\tilde X_{t*} e_i| \<N_t,\partial_{r_{p_0}}\>,$$  so
\be
 |\tilde X_{t*} e_i|  \le \frac{\s(\r)}{\<N,\partial_{r_{p_0}}\>} < \frac{\s(\rho_t + \ral \ln 2)}{\ral\ \ta(\rho_t)} \leq  \frac{\s(\psi(V_0)+ \ral \ln 2)}{\ral\ \ta(\xi(\psi(V_0)))}, 
\ee
where we have used Theorem \ref{TBGRM} for the second inequality and \eqref{brho} for the third one. Moreover, it follows from \eqref{Xtei}  that  $|e_i(\r)| \le  |\tilde X_{t*} e_i| $, thus both the first derivatives of $\tilde X_t$ and $\r$ are bounded independently of $t$.

On the other hand, it is clear from the expression \eqref{Xtilde} for $\tilde X_t$ that all the higher order derivatives of $\tilde X_t$ are bounded if an only if the corresponding derivatives of $\r$ are bounded. In order to see that such derivatives of $\r$ are bounded, first we compute the components $\alpha_{ij}$ of the second fundamental form of $\varphi_t(M_t)$ using the parametrization \eqref{Xtilde}, that is,  $\alpha_{ij} = \alpha(\tilde X_{t*} e_i, \tilde X_{t*} e_j)$.
We shall write $\alpha_{ij}$ in terms of $\r$ and its derivatives.

 If $\xi$ is a vector normal to $\varphi_t(M_t)$ satisfying $\<\xi, \partial_ {r_{p_0}}\> = s_{\lambda}(r_{p_0})$, we have 
 $$ 0 =\<\xi, \tilde X_{t*} e_i\> \circ \tilde X_t = e_i(\r) s_\lambda(\r) +  s_\lambda(\r) \,\<\tau_s e_i, \xi\>,$$
and then (without explicit writing of the suitable  compositions with the map $\tilde X_t$)
$$\<\tau_s e_i, \xi\> = - e_i(\r), \quad \text{so} \quad \xi =  s_\lambda(\r) \partial_{r_{p_{0}}} - \sum_{i =1}^n e_i (\r) \tau_s e_i.$$
Consequently, the outward  unit normal vector $N_t$ to $\varphi_t(M_t)$ can be written as
\begin{align} \label{normal}
&N_t= \frac1{ |\xi|} \left(s_\lambda(\r) \partial_{r_{p_{0}}}   - \sum_{i =1}^n e_i (\r) \tau_s e_i\right), \text{ with } \\
&|\xi| =\sqrt{s_{\lambda}^2(\r) + |\!\gras \! \r|^2}.\nonumber
\end{align}

To compute the components $\alpha_{ij}$, we use on $\eme$ the spherical coordinates $\gamma: \re^+\times S^n \flecha \eme$ defined by $\gamma(s,u) = \exp_{p_0} s\, u$. In these coordinates, for a local orthonormal frame $\{E_0= \partial_{r_{p_0}}, E_1=\tau_s e_1, ..., E_n=\tau_s e_n\}$ of $\eme$ and its dual frame $\{dr_{p_0}, \theta^1, ..., \theta^n\}$, we have $\gamma^* d {r_{p_0}} = ds$, $\gamma^*\theta^i = \s e^i$, being $\{e^1, ... , e^n\}$ the dual frame of $\{e_1, ..., e_n\}$. Let us denote by $\nabla_S$ and $g_S$  the standard covariant derivative and  metric of  $S^n$, respectively.

The Cartan connection $1$-forms $\omega_0^j$, $\omega_i^j$ of $\nablaa$    satisfy $\omega_0^j = - \co \theta^j$ and $\gamma^*\omega_i^j = {\mf{s}}_i^j$, where ${\mf s}_i^j$ are the connection forms of $\nabla_S$ in the frame $\{e_1, ..., e_n\}$.  Using these facts and \eqref{normal},  after a standard computation, we reach
\begin{align}   \label{2FF_rad}
\alpha_{ij} &= - \<\nablaa_{\tilde X_{t*}e_i} {\tilde X_{t*}e_j}, N_t\> \\
&= -\frac1{|\xi|}
\left(s_\lambda(\r) \nabla_S^2 \r - s_\lambda^2(\r) c_\lambda(\r) g_S -  2 c_\lambda(\r) d\r \otimes d\r\right)(e_i,e_j). \nonumber 
\end{align} 

Since each $\varphi_t$ is an isometry of the ambient space, the second fundamental forms of $M_t$ and $\varphi_t(M_t)$ coincide. Then, by \eqref{Lmbound}, $\alpha$ and all their derivatives are uniformly bounded and, by \eqref{2FF_rad} and the fact that $\r$ and its first order derivatives are uniformly bounded, we have that all the derivatives of $\r$ are uniformly bounded. Thus, by the relation \eqref{Xtilde}, all the derivatives of $\tilde X_t$ are also uniformly bounded. 

We are now in conditions to  apply Arzelˆ-Ascoli Theorem to conclude the existence of  sequences of maps $\tilde X_{t_i}$ and $\tilde r_{t_i}$ satisfying \eqref{Xtilde} which $C^\infty$-converge to  smooth  maps $\tilde X_\infty: S^n \flecha \eme$ and  $\tilde r_\infty : S^n\flecha \re^+$ satisfying $\tilde X_\infty(u) = \exp_{p_0} \tilde r_\infty(u) u$. The last equality implies that $\tilde X_\infty$ is an immersion and, since the convergence is smooth and all the hypersurfaces $\tilde X_t(S^n)$ are $h$-convex, we can assure that $\mc S = \tilde X_\infty(S^n)$ is $h$-convex. Using  Remark (iii) in the Introduction, we say that $\varphi_{t_i}(M_{t_i})$ converges to $\mc S$ as $t_i \to \infty$.
 
 Finally, by Step 1, $\mc S$  must be a compact embedded hypersurface in $\eme$ of constant mean curvature, that is, a geodesic sphere of $\eme$.
This finishes the proof of Proposition \ref{SeqConv}.  
\end{demo}

\section{Exponential convergence}\label{ExCo}

In order to complete the proof of statement (c) in Theorem \ref{math}, 
our next goal is to show that the $M_t$'s converge to some limit $M_\infty$ exponentially. First, let us fix an instant $t_k\in [0, \infty[$. We can parametrice $M_{t}$, with $t \geq t_k$, by 
\be \label{param1}
X_t(x) = \exp_{p_{t_k}} r(t, u(t,x)) u(t,x),
\ee
  where $u(t,x) = \ds\frac{\exp_{p_{t_k}}^{-1}X_t(x)}{r_{p_{t_k}}(X_t(x))}  \text{ and }  r(t, u(t,x)) = r_{p_{t_k}}(X_t(x)).$
At least for $t$ near to $t_k$, we have $p_{t_k} \in \Omega_t$, and so the map $u_t: M \flecha S^n\subset T_{p_{t_k}} \eme$ defined by $u_t(x)=u(t,x)$ is a diffeomorphism. 

Observe that the map
\be
\overline X_t(x) = \exp_{p_{t_k}} r(t, u(t_k,x))\  u(t_k,x)
\ee
is another parametrization of $M_t$. In fact, writing $\phi_t= u_{t_k}^{-1} \circ u_t :M \flecha M$, we see that $\overline X_t \circ \phi_t = X_t$, i.e., the motions $X_t$ and $\overline X_t$ differ only by a tangential diffeomorphism $\phi_t$. Moreover, 
\be 
 \parcial{X_t}{t} = \parcial{\overline X_t }{t} \circ \phi_t  + \overline X_{t*}\parcial{\phi_t}{t}.
\ee
Therefore, if $X_t$ is  a solution of  \eqref{vpmf}, $\overline X_t$ satisfies the equation
\be\label{evol.oX}
\<\parcial{\overline X_t}{t} , N_t\> = (\oH_t - H_t).
\ee
Conversely,  it is well known (see, for instance, \cite{Eck}) that \eqref{evol.oX} is equivalent to \eqref{vpmf} (by tangential diffeomorphisms).

With the aim of applying some methods from \cite{EsSi98}, it is convenient to write \eqref{evol.oX} as an equation for the function $r(t, \cdot\,)$. Previously, to simplify the notation, we shall compose with the diffeomorphism $u_{t_k}^{-1}$ in order to consider $\overline X_t$ as a map from $S^n$ (instead of $M$) into $\eme$, i.e., 
\be\label{oXtex}
\overline X_t(u) = \exp_{p_{t_k}} r(t,u)\ u \quad   \text { for every } u\in S^n.
\ee 

For any local orthonormal frame $\{e_i\}$ of $S^n$, a basis of the tangent space to $M_t$ is given  by
$\{\tilde{e}_i = \overline X_{t \ast} e_i \}$. In this basis,  the outward  unit normal vector $N_t$ to $M_t$ and the second fundamental form $\alpha_t$ are given by the expressions \eqref{normal} and \eqref{2FF_rad}, respectively (with the obvious change of $\r$ by $r(t,\cdot)$).

From \eqref{evol.oX}, \eqref{oXtex} and \eqref{normal}, we obtain
\be 
\parcial{r}{t} = s_\lambda^{-1}(r)  (\oH_t - H_t) \sqrt{s_\lambda^2(r) + |\! \gras \! r|^2} \label{rel_r_H}
\ee

On the other hand, the components of the metric $g_t$ in the basis $\{\tilde e_i\}$ are
$$g_{ij} = e_i(r) e_j(r) + s_\lambda^2(r) \delta_{ij}$$
From this, using an elementary algebraic result, we can express the components of the inverse metric as
\be \label{inv_met}
g^{ij} = \frac1{s_\lambda^2(r)} \left( \delta^{ij} - \frac1{|\xi|^2} e_i(r) e_j(r)\right),
\ee

Then, joining \eqref{inv_met} and \eqref{2FF_rad}, we get
\begin{align} \label{H_rad}
H_t  = & - \frac{s_\lambda^{-1}(r)}{|\xi|} \left(\Delta_S r - \frac1{|\xi|^2} \nabla_S^2r(\grad_{S^n}r,\grad_{S^n}r) \right)  \\
&+ \frac{c_\lambda(r)}{|\xi|} \left(n + \frac{|\! \gras \! r|^2}{|\xi|^2}\right) \nonumber
\end{align}

Finally, substituting \eqref{H_rad} in \eqref{rel_r_H}, we can write
\begin{align} 
\parcial{r}{t} = &  s_\lambda^{-2}(r) \left(\Delta_S r - \frac1{|\xi|^2} \nabla_S^2r(\grad_{S^n}r,\grad_{S^n}r) \right) \label{ec_G}\\
& - co_\lambda(r) \left(n + \frac{|\! \gras\! r|^2}{|\xi|^2}\right) + s_\lambda^{-1}(r) \oH_t |\xi|. \nonumber 
\end{align}

Observe that equation \eqref{ec_G} coincides with equation (2.1) in \cite{EsSi98} when we change $\s(r)$ by $r$ and $\c(r)$ by $1$. Therefore, \eqref{ec_G} satisfies all the conditions which allow to apply (vii) in \cite{EsSi98}, and conclude

\begin{prop}\label{rKe}
Given $m\in \mathbb{N}$ and a constant $s > 0$, there exists $\omega>0$ and a neighborhood $V $ of $s$ in $h^{1 + \beta}(S^n)$ such that for each initial condition $r_0 \in V$

(a) The solution $r(t, \cdot \,)$ of \eqref{ec_G} satisfying $r(0,\cdot \,) = r_0(\cdot \,)$ exists on $[0,\infty[$, and

(b) there exist $c=c(m,\omega) > 0$, $T=T(m,\omega)>0$,   a unique function $\tilde \rho$ (in a space of functions $ \mc{M}^c$ called center manifold)  and $K=K(r_0, c, \tilde\rho)$ such that
$$\|r(t,\cdot \,) - \tilde \rho(\cdot \,)\|_{C^m} \leq K\  e^{-\omega t} $$
for $t > T$.
\end{prop}

Now we are in position to finish the proof of the exponential convergence of $M_t$ to a geodesic sphere. Indeed, 

\medskip
\begin{demo} {\bf of (c), Theorem 1} \
Let us apply Proposition \ref{SeqConv} to take $t_k$ big enough so that $\varphi_{t_k}(M_{t_k})$  is near to the limit  geodesic sphere $\mc S$ of radius $\mf r$. As $\varphi_{t_k}$ is an isometry of $M^{n+1}_\lambda$, we have that $M_{t_k}$ is close to the geodesic sphere $\varphi_{t_k}^{-1}(\mc{S})$ of radius $\mf r$. Thus, using spherical geodesic coordinates, we can assure that  $r_{p_{t_k}}$ belongs to a small neighborhood $V$ of the constant function $\mf r$ in $h^{1+\beta}(S^n)$.

So, applying Proposition \ref{rKe} with initial condition $r_{p_{t_k}}$, we can conclude that  the solution $r(t, \cdot \,)$ of \eqref{ec_G}  starting at $r_{p_{t_k}}$ is defined on $[0, \infty[$ and converges exponentially to a unique function $\tilde \rho$.  This implies that $\overline X_t(u) = \exp_{p_{t_k}} r(t,u) u$ solves \eqref{evol.oX} and converges exponentially to $u\mapsto \exp_{p_{t_k}} \tilde\rho(u) u$. Therefore,  the reparametrization $X_t$ of $\overline X_t$ given by  \eqref{param1} has the same convergence properties; in addition, it is a solution of \eqref{vpmf} starting at $r_{p_{t_k}}$, and, by uniqueness, $X_t$ coincides on $[t_k, \infty[$ with the solution of \eqref{vpmf} given by part (b) of Theorem \ref{math}.

On the other hand, Step 1 in the proof of Proposition \ref{SeqConv} says that the mean curvature $H_t$ of the hypersurfaces $X_t(M)$ tends to a constant value as $t\to \infty$. In conclusion,
 the only possibility is that $\exp_{p_{t_k}} \tilde\rho(u) u$ represents a geodesic sphere in $M^{n+1}_\lambda$ and, by the volume-preserving properties of the flow, such sphere has to enclose the same volume as the initial condition $X_0(M)$.
\end{demo}

\section{A result for certain non-necessarily $h$-convex initial data}

A remarkable fact is that in the last section we have not used all the strength of the results on the existence and exponential attractivity of the center manifold $\mc{M}^c$. It is precisely this additional power which allows us to extend the claims about long time existence and convergence of Theorem \ref{math} to certain non $h$-convex initial data; in particular, those \emph{sufficiently close} to a geodesic sphere of $M^{n+1}_\lambda$.

Notice that, if we begin the flow with a non $h$-convex hypersurface $M_0 \subset M^{n+1}_\lambda$, we cannot establish the convergence of a sequence $\{M_{t_i}\}$ up to isometries (like in Proposition \ref{SeqConv}), because the properties of $h$-convexity are strongly used to find  $t$-independent bounds for the second fundamental form (together with all its derivatives) of the hypersurfaces $M_t$ evolving under the flow, and recall that these bounds are the key to prove \eqref{HtooH}.

In spite of this, it is not difficult to overcome the absence of $h$-convexity since we are just in the same situation as in \cite{EsSi98} (see also \cite{E-S 2} for a full understanding of \cite{EsSi98}). The 
 only point which need to be checked again in our particular situation is  the equality between the center manifold $\mc{M}^c$ (cf. \cite{EsSi98} for its definition) and the equilibria $\mc M$ of \eqref{ec_G} in some small neighborhood, as it is obvious that $\mc M$ is different in equations \eqref{ec_G} and \cite{EsSi98} (2.1). Next we are going to check this identity.

\begin{prop} \label{M=Mc} 
Let $\mc{S}$ be a geodesic sphere of $M_\lambda^{n +1}$ of radius $r_{_{\! \! \mc{S}}}$ and center $p_{_{\! \mc{S}}}$. There is a neighborhood $\mc{O}$ of $r_{_{\! \! \mc{S}}}$ in  which $\mc{M}$ coincides with an open set of the local center manifold $\mc{M}^c$ for the equation \eqref{ec_G}. 
\end{prop}

\begin{demo} 
Let us begin by observing that the construction of $\mc{M}^c$ as a center manifold for \eqref{ec_G} is identical to that for the equation (2.1) in \cite{EsSi98}. Therefore, $\mc{M}^c$ is a $(n+2)$-dimensional manifold tangent to $ \{1\}\oplus \mc{H}_1$, where $ \{1\}$ denotes the space of constant functions on $S^n$ and $\mc{H}_1$ is the space of eigenfunctions corresponding to the first nonzero eigenvalue of $\Delta_{S}$. 

On the other hand, $\mc M$ is the space of functions $\rho:S^n\flecha \re^+$ such that $\exp_{p_{_\mc{S}}}\rho(u) u$ parametrices a constant mean curvature hypersurface of $\eme$, that is, a geodesic sphere of $\eme$. Then, in a small neighborhood  of $r_{_{\! \! \mc{S}}}$
\be\label{Ue}
\mathbb{U}_\eps=\{ r \in C^\infty(S^n)\, : \, ||r - r_{_{\! \! \mc{S}}}||_{h^{1+\beta}} < \eps\},
\ee 
 $\mc M$ can be parametrized by $(z_0,z)\in \re^{n+2} \equiv \re \oplus T_{p_{_{\! \mc{S}}}}\eme$, being $z_0 + r_{_{\! \! \mc{S}}}$ the radius of a geodesic sphere $S_z$ and $z=(z_1, ..., z_{n+1})$ the normal coordinates (centered at $p_{_{\! \mc{S}}}$) of its center. 

The function $\rho_z:S^n\flecha \re^+$ which represents the geodesic sphere $S_z$ has to satisfy $r_{_{\! \! \mc{S}}} + z_0 = \dist(\exp_{p_{_{\! \mc{S}}}}\rho_z(u) u, \exp_{p_{_{\! \mc{S}}}}z) $. Using hyperbolic trigonometry, we can write this equality under the form
\be\label{ec_esfM}
\c(r_{_{\! \! \mc{S}}} + z_0) = \c(\rho_z(u)) \c(|z|)  + \lambda \, \s(\rho_z(u)) \frac{\s(|z|)}{|z|} \<u,z\>.
\ee
Thus, by implicit differentiation of \eqref{ec_esfM}, we obtain, at $(z_0,z)=(0,0) \in \re^{n+2}$,
\be\label{drz}
\left.\parcial{\rho_z}{z_0}\right|_{(0,0)}(u) = 1, \quad \left.\parcial{\rho_z}{z_i}\right|_{(0,0)}(u) =  u_i.
\ee

Since $\{1, u_1, \ldots , u_{n +1}\}$ is a basis of $\{1\} \oplus \mc{H}_1$, the differential of the function $\rho: (z_0,z)\mapsto\rho_z$ at $(0, 0)$ is an isomorphism between $ \re^{n+2}$ and  $\{1\} \oplus  \mc{H}_1$. From now on, the equality $\mc{M} = \mc{M}^c$ in a neighborhood $\mc{O}$ of $r_{_{\! \! \mc{S}}}$ follows like in (vi) of \cite{EsSi98}. 
\end{demo}

\begin{nota} \label{ntma2}
Arguing as in (vii) of \cite{EsSi98}, we can conclude that, given $m \in \mathbb{N}$, there exists $\omega >0$ and $\eps > 0$ such that for each initial condition $r_0 \in \mathbb{U}_\eps$, where $\mathbb{U}_\eps$ is defined by \eqref{Ue}, it follows the statement (a) in Proposition \ref{rKe}. Moreover, as in part (b) of the same proposition, we can find constants   $c=c(m,\omega) > 0$, $T=T(m,\omega)>0$, and a unique $\tilde{\rho} \in \mc{O} \cap \mc{M}^c$, depending only on $r_0$, and satisfying
$$\|r(t,\cdot \,) - \tilde \rho(\cdot \,)\|_{C^m} \leq K(r_0, c,\tilde\rho)\  e^{-\omega t}  \qquad \text{for all } \quad t > T.$$

But applying Proposition \ref{M=Mc}, we know that $\mc{M}^c \cap \mc{O} = \mc{M}$; so we can find a unique $(z_0, z) \in \rho^{-1}(\mc{O})$ such that $\exp_{p_{_{\! \mc{S}}}} \tilde \rho (u) u$ represents a geodesic sphere of $M^{n+1}_\lambda$ with center $z$ and radius $r_{_{\! \! \mc{S}}} + z_0$.
\end{nota}

\medskip

Thanks to Proposition \ref{M=Mc}, we are in position to prove Theorem \ref{corm}.

\medskip

\begin{demo} {\bf of Theorem \ref{corm}} \ 

Let $\mc{S} = \partial B(p_{_{\! \mc{S}}}, r_{_{\! \! \mc{S}}})$ be a geodesic sphere in $M^{n+1}_\lambda$ and $m\in \mathbb{N}$. By Remark \ref{ntma2}, we can find an $\eps > 0$ and a neighborhood $\mathbb{U}_\eps$, defined as in \eqref{Ue}, satisfying the property detailed in Remark \ref{ntma2}.

Now consider any arbitrary embedding $X: M \fle M_\lambda^{n+1}$ which $h^{1 + \beta}$-distance to $\mc{S}$ is less than $\eps$; in other words, taking spherical coordinates centered at $p_{_{\! \mc{S}}}$, the radial distance $r(\cdot\,) = dist(p_{_{\! \mc{S}}}, X(\cdot \,))$ of $X(M)$ belongs to the neighborhood $\mathbb{U}_\eps$.

Therefore, Remark \ref{ntma2} assures that the solution $r_t(\cdot\,)$ of \eqref{ec_G} starting at $r$ exists on $[0, \infty[$ and $X_t(u) = \exp_{p_{_{\! \mc{S}}}} r_t(u) u$ converges, as $t \to \infty$, to a geodesic sphere in $M_\lambda^{n+1}$.
Finally, noticing that $X_t$ is a solution of \eqref{evol.oX}, Theorem \ref{corm} follows.
\end{demo}

\medskip
{\bf Acknowledgments.}
We want to thank F.J. Carreras and M. RitorŽ for their help at different points of this paper.

Both authors have been   partially supported by the \uppercase{DGI(E}spa–a) and \uppercase{FEDER P}roject \uppercase{MTM}  \uppercase{N}o {\it 2004-06015-\uppercase{C}02-01}. The first author has been supported by the \uppercase{B}eca del  \uppercase{P}rograma \uppercase{N}acional de \uppercase{F}ormaci—n del \uppercase{P}rofesorado \uppercase{U}niversitario ref: \uppercase{AP}2003-3344.

{\footnotesize

\bibliographystyle{alpha}

}

\vskip1truecm

{\small 

\begin{tabular}{ c c }
Address & Universidad de Valencia.\\
      \       & 46100-Burjassot (Valencia) Spain \\
      \       & email: esther.cabezas@uv.es\\
      \       & and\\
      \       & email: miquel@uv.es 
\end{tabular}

}

\end{document}